\documentclass[russian]{article}
\usepackage[T1,T2A]{fontenc}
\usepackage[cp1251]{inputenc}
\usepackage{babel}
\usepackage{mathtools}
\usepackage{amsmath}
\usepackage{amssymb}
\usepackage{esint}
\usepackage[unicode=true] {hyperref}



\begin{document}

\title{The sojourn time problem for a $p$-adic random walk and its applications to the spectral diffusion of proteins}
\author{A.\,Kh.~Bikulov \\
 \textit{Institute of Chemical Physics, } \\
 \textit{Kosygina Street 4, 117734 Moscow, Russia} \\
 e-mail:\:\texttt{beecul@mail.ru} \\
 and \\
 A.\,P.~Zubarev \\
 \textit{ Physics Department, Samara State Aerospace University, } \\
 \textit{ Moskovskoe shosse 34, 443123, Samara, Russia} \\
 \textit{Physics and Chemistry Department, } \\
 \textit{Samara State University of Railway Transport,} \\
 \textit{Perviy Bezimyaniy pereulok 18, 443066, Samara, Russia} \\
 e-mail:\:\texttt{apzubarev@mail.ru} }
\maketitle
\begin{abstract}
We consider the problem of the distribution of the sojourn time in
a compact set $\mathbb{Z}_{p}$ in the case of a $p$-adic random
walk. We rely on the results of our previous studies of the
distribution of the first return time for a $p$-adic random walk
and the results of Takacs on the study of the sojourn time problem
for a wide class of random processes. For a $p$-adic random walk
we find the mean sojourn time of the trajectory in
$\mathbb{Z}_{p}$ and the asymptotics as $t\rightarrow\infty$ of
arbitrary moments of the distribution of the sojourn time in
$\mathbb{Z}_{p}$. We also discuss some possible applications of
our results to the modeling of relaxation processes related to the
conformational dynamics of protein.

\vspace{5mm}
{\bf Keywords:} $p$-adic mathematical physics,  $p$-adic random walk, sojourn time problem, $p$-adic models of conformational dynamics
\end{abstract}

\section{Introduction}

For a wide class of random processes, the problem of finding the
sojourn time of the trajectory of a process in an arbitrary subset
was solved by Takacs \cite{Takacs_1,Takacs_2,Takacs_3} and is
currently widely used in various applications of the theory of
probability and random processes (see, for example,
\cite{BP,Io,Gr}). In the present study, we apply these results to
a random walk on the field of $p$-adic numbers $\mathbb{Q}_{p}$ (a
$p$-adic random walk).

The analysis of $p$-adic random processes attracts keen interest
in relation to the study of processes and phenomena that have an
explicit or hidden ultrametric structure in various fields of
physics, biology, and informatics: spin glasses, proteins,
cognitive processes, genetic code, biological evolution,
classification, taxonomy, and other fields (see, for example,
\cite{RTV,VVZ,ALL,ALL_1,DKKM}). It is worth special mention that a
$p$-adic random walk provides an adequate mechanism to describe
the dynamics of conformational rearrangement of protein molecules
in the native state. The problem, considered in the present paper,
of the sojourn time of the trajectory of a random walk in a given
domain is directly related to such processes. It is known that
ultrametric modeling of the conformational dynamics of protein is
based on the description of the dynamics of a system as a random
process in the ultrametric state space of the system, where
probabilistic transitions between states are defined by the
ultrametric distance between them. From the mathematical
viewpoint, such a dynamics is described by the equation of
$p$-adic random walk with the Vladimirov operator \cite{VVZ} with
possible reaction terms. Similar models admitting exact analytic
solutions were considered in
\cite{ABK_1999,ABKO_2002,ABO_2003,ABO_2004,AB_2008,ABZ_2009,ABZ_2011,ABZ_2013,BZ_2021}
as models of ultrametric random walk on energy landscapes of
protein molecules, models of reaction of myoglobin rebinding to
small ligands, models of molecular nano-machine prototypes, etc.
In paper \cite{AB_2008} devoted to the description of experiments
on the spectral diffusion in proteins, the variation in the
absorption frequency of a chromophore marker bound to the active center
of protein is associated with the average number of returns of the
protein trajectory to a given domain in the conformational space
that is directly related to the sojourn time of the trajectory in
the given domain of conformations. Therefore, in the general case,
the problem of the sojourn time of protein in a given
conformational domain is directly related to the description of
processes in protein that can occur only if the protein is in a
certain conformational subspace.

In this paper, following our previous studies, we consider a
$p$-adic random walk as a Markov process
$\xi\left(t,\omega\right):R_{+}^{1}\times\Omega\to\mathbb{Q}_{p}$
($\Omega$ is a set of elementary events with a sigma algebra
$\Sigma$ and a probability measure ${\rm P}$) with the transition
function $f(y,t|x,0)\equiv f(y-x,t)$. Below we will omit the
argument $\omega$ and write
$\xi\left(t\right)\equiv\xi\left(t,\omega\right)$. The function
$f\left(x,t\right)$ is the fundamental solution of the equation of
$p$-adic random walk with the Vladimirov operator
\begin{equation}
\frac{\partial}{\partial t}f\left(x,t\right)=-\frac{1}{\Gamma_{p}\left(-\alpha\right)}\intop_{\mathbb{Q}_{p}}d_{p}y\frac{f\left(y,t\right)-f\left(x,t\right)}{\left|y-x\right|_{p}^{\alpha+1}},\label{UDE}
\end{equation}
where $\Gamma_{p}(-\alpha)=\dfrac{1-p^{-\alpha-1}}{1-p^{\alpha}}$
is the $p$-adic gamma function and the integration is performed
with respect to the Haar measure on the field $\mathbb{Q}_{p}$
\cite{VVZ}. Since the solution $f_{r,a}\left(x,t\right)$ of
equation (\ref{UDE}) with the initial condition in an arbitrary
ball $f_{r,a}\left(x,0\right)=\Omega(|x-a|_{p}p^{-r})\equiv\left\{
\begin{array}{l}
1,{\rm \;}|x-a|_{p}\le p^{r},\\
0,{\rm \;}|x-a|_{p}>p^{r}.
\end{array}\right.$ is scale- and translation-invariant, i.e., $f_{r,a}\left(x,t\right)=f_{0,0}\left(p^{-r}\left(x-a\right),p^{-\alpha r}t\right),$
we will assume without loss of generality that the initial
distribution of the random process $\xi\left(t\right)$ has the
form

\begin{equation}
\varphi\left(x,0\right)=\Omega(|x|_{p})\label{IC}
\end{equation}
In this case, the densities of the distribution function
$\varphi(x,t)$ and the transition function $f(y-x,t)$ of the
process $\xi\left(t\right)$ are, respectively,

\[
\varphi(x,t)=\intop_{\mathbb{Q}_{p}}\Omega\left(\left|k\right|_{p}\right)\exp\left(-\left|k\right|_{p}^{\alpha}t\right)\chi\left(-kx\right)d_{p}k,
\]

\[
f(y-x,t)=\intop_{\mathbb{Q}_{p}}\exp\left(-\left|k\right|_{p}^{\alpha}t\right)\chi\left(-k\left(y-x\right)\right)d_{p}k,
\]
where $\chi\left(x\right)$ is an additive character of the field
$\mathbb{Q}_{p}$ \cite{VVZ}.

For the process $\xi\left(t\right)$, consider random processes
$I_{\mathbb{Z}_{p}}\left(t\right)$ and
$I_{\mathbb{Q}_{p}\setminus\mathbb{Z}_{p}}\left(t\right)$, defined
as
\[
I_{\mathbb{Z}_{p}}\left(t\right)=\left\{ \begin{array}{c}
1,\;\xi\left(t\right)\in\mathbb{Z}_{p}\\
0,\;\xi\left(t\right)\in\mathbb{Q}_{p}\setminus\mathbb{Z}_{p}
\end{array}\right.,\;I_{\mathbb{Q}_{p}\setminus\mathbb{Z}_{p}}\left(t\right)=\left\{ \begin{array}{c}
1,\;\xi\left(t\right)\in\mathbb{Q}_{p}\setminus\mathbb{Z}_{p}\\
0,\;\xi\left(t\right)\in\mathbb{Z}_{p}
\end{array}\right.,
\]
as well as random processes
$\theta_{\mathbb{Z}_{p}}\left(t\right)$ and
$\theta_{\mathbb{Q}_{p}\setminus\mathbb{Z}_{p}}\left(t\right)$:
\begin{equation}
\theta_{\mathbb{Z}_{p}}\left(t\right)=\intop_{0}^{t}I_{\mathbb{Z}_{p}}\left(t'\right)dt,\label{th_Z_p}
\end{equation}
\begin{equation}
\theta_{\mathbb{Q}_{p}\setminus\mathbb{Z}_{p}}\left(t\right)=\intop_{0}^{t}I_{\mathbb{Q}_{p}\setminus\mathbb{Z}_{p}}\left(t'\right)dt.\label{th_Q_p-Z_p}
\end{equation}
We will call the processes (\ref{th_Z_p}) and (\ref{th_Q_p-Z_p})
the sojourn time of the trajectory in $\mathbb{Z}_{p}$ $\text{and
}\mathbb{Q}_{p}\setminus\mathbb{Z}_{p}$ in the time interval from
$0$ to $t$. The goal of the present work is the study of the
processes (\ref{th_Z_p}) and (\ref{th_Q_p-Z_p}). In solving this
problem, we use the results of our paper \cite{ABZ_2009} devoted
to the solution of the problem of first return to the domain
$\mathbb{Z}_{p}$, as well as the results of
\cite{Takacs_1,Takacs_2,Takacs_3} devoted to the problem of the
sojourn time of a trajectory of a wide class of random processes
in an arbitrary subset.

The paper is organized as follows. In Section 2, we formulate the
main results of solving the first return problem that we obtained
in \cite{ABZ_2009} and use in the present study. In Section 3, we
adapt the main theorem on the sojourn time, which was proved in
\cite{Takacs_1,Takacs_2,Takacs_3}, to the case of a $p$-adic
random walk. In Section  4, for a $p$-adic random walk we find the
mean sojourn time of the trajectory in $\mathbb{Z}_{p}$ and
asymptotics as $t\rightarrow\infty$ of the moments of the
distribution of the sojourn time in $\mathbb{Z}_{p}$. In the
concluding section, we discuss some possible applications of our
results to modeling relaxation processes related to the
conformational dynamics of protein.

\section{First return to the domain $\mathbb{Z}_{p}$ }

Our solution of the problem of the sojourn time of the trajectory
of a $p$-adic random walk in a given domain is based on the
results we obtained when solving the problem of the distribution
of the first return time to the domain $\mathbb{Z}_{p}$ for the
process $\xi\left(t\right)$ in \cite{ABZ_2009}. In this section,
we present the main results of that paper, which will be needed in
what follows.

The first return time of the trajectory of a random process
$\xi\left(t\right)$ to the domain $\mathbb{Z}_{p}$ is a random
variable $\tau_{Z_{p}}\in\mathbb{R}_{+}$ defined as
\[
\tau_{Z_{p}}=\inf\left\{ t>0:\:\exists t':\;0<t'<t,\:\xi\left(t'\right)\in\mathbb{Q}_{p}\setminus\mathbb{Z}_{p},\:\xi\left(t\right)\in\mathbb{Z}_{p}\right\} .
\]
In general, the first return time $\tau_{Z_{p}}$ is a random
variable in an extended sense and can take values on extended
nonnegative number axis
$\mathbb{\bar{R}}_{+}=\mathbb{R}_{+}\cup\left\{ +\infty\right\} $.
A process $\xi\left(t\right)$ is said to be recurrent if
$\tau_{Z_{p}}$ exists and its values lie in $\mathbb{R}_{+}$ with
probability $1$. If $\tau_{Z_{p}}$ exists and takes values in
$\mathbb{R}_{+}$ with probability $q\text{<1}$, then the process
$\xi\left(t\right)$ is said to be nonrecurrent. For a recurrent
process, the distribution density $f\left(t\right)$ of
$\tau_{Z_{p}}$ is normalized to $1$. For a nonrecurrent process,
$\intop_{0}^{\infty}f\left(t\right)dt=q$.

\noindent \textbf{Theorem 1}. The distribution density
$f\left(t\right)$ of the first return time $\tau_{Z_{p}}$
satisfies the inhomogeneous Volterra equation
\begin{equation}
v(t)=\intop_{0}^{t}v(t-\tau)f(\tau)d\tau+f(t)\label{g}
\end{equation}
where
\[
v(t)=-\frac{1}{\Gamma_{p}(-\alpha)}\intop_{\mathbb{Q}_{p}\setminus\mathbb{Z}_{p}}\frac{\varphi(x,t)}{\left|x\right|_{_{p}}^{\alpha+1}}dx
\]
is the probability of transition in unit time from
$\mathbb{Q}_{p}\setminus\mathbb{Z}_{p}$ to $\mathbb{Z}_{p}$ at
time $t$.

\noindent The proof of this theorem is given in
\cite{ABZ_2009}.\textbf{ }

\noindent Denote by
\[
J\left(t\right)=\intop_{\mathbb{Z}_{p}}\varphi\left(x,t\right)d_{p}x
\]

\noindent the probability to find the trajectory in
$\mathbb{Z}_{p}$ at time $t$ (we will call it the survival
probability). Then (\ref{UDE}) implies
\begin{equation}
\dfrac{dJ\left(t\right)}{dt}=v(t)-B_{\alpha}J\left(t\right)\label{dJ}
\end{equation}
where
\[
B_{\alpha}=-\dfrac{1}{\Gamma_{p}\left(-\alpha\right)}\intop_{\mathbb{Q}_{p}\setminus\mathbb{Z}_{p}}\dfrac{1}{\left|x\right|_{p}^{\alpha+1}}dx=\dfrac{1-p^{-1}}{1-p^{-\alpha-1}}.
\]
Denoting the Laplace transforms of the functions $v(t)$, $f(t)$,
and $J\left(t\right)$, respectively, by $\hat{v}(s)$,
$\hat{f}(s)$, and $\hat{J}\left(s\right)$, from (\ref{g}) we
obtain
\begin{equation}
\hat{f}(s)=\frac{\hat{v}(s)}{1+\hat{v}(s)}.\label{F}
\end{equation}

\noindent In terms of Laplace transforms, (\ref{dJ}) has the form

\begin{equation}
\hat{v}(s)=(B_{\alpha}+s)J(s)-1,\label{v_L}
\end{equation}
where $\hat{J}\left(s\right)$ is the Laplace transform of
$J\left(t\right)$:
\[
\hat{J}\left(s\right)=\left(1-p^{-1}\right)\sum_{n=0}^{\infty}\dfrac{p^{-n}}{s+p^{-\alpha n}}.
\]
It follows from the representation (\ref{v_L}) that
\begin{equation}
\hat{f}\left(s\right)=1-\frac{1}{\left(B_{\alpha}+s\right)\hat{J}\left(s\right)}.\label{F_1}
\end{equation}

Asymptotic estimates for the density of the distribution function
$f(t)$ of the first return time are given by the following
theorem, which was proved in \cite{ABZ_2009}:

\noindent \textbf{Theorem 2. } The function $f(t)$ has the
following asymptotic estimates as $t\to\infty$:
\begin{eqnarray*}
A_{1}t^{-\frac{2\alpha-1}{\alpha}}\left(1+o(1)\right) & \le & f\left(t\right)\le B_{1}t^{-\frac{2\alpha-1}{\alpha}}\left(1+o\left(1\right)\right)\;{\rm for}\;\alpha>1,\\
A_{2}t^{-\frac{1}{\alpha}}\left(1+o(1)\right) & \le & f\left(t\right)\le B_{2}t^{-\frac{1}{\alpha}}\left(1+o\left(1\right)\right)\;\mbox{{\rm for}}\;\alpha<1,\\
A_{3}\frac{t^{-1}}{(\log t)^{2}}\left(1+o(1)\right) & \le & f\left(t\right)\le B_{3}\frac{t^{-1}}{(\log t)^{2}}\left(1+o\left(1\right)\right)\;\mbox{{\rm for}}\;\alpha=1,
\end{eqnarray*}
where $o\left(1\right)\to0$ as $t\to\infty$ and $A_{i}$ and
$B_{i}$, $i=1,2,3$ are some functions of $\alpha$ and $p$.

Note that the distribution function $f\left(t\right)$ of the first
return time has the following property. For $\alpha\ge1$, we have
$\intop_{0}^{\infty}f\left(t\right)dt=\hat{f}\left(0\right)=1$,
which means that the trajectory returns with probability 1 to the
domain $\mathbb{Z}_{p}$ in finite time. Hence, for $\alpha\ge1$,
the random walk is recurrent. For $0<\alpha<1$, we have
\[
\intop_{0}^{\infty}f\left(t\right)dt=\hat{f}\left(0\right)=\dfrac{p}{p^{\alpha}}\left(\dfrac{p^{\alpha}-1}{p-1}\right)^{2}\equiv C_{\alpha}<1,
\]
which means that, with probability $1-C_{\alpha}$, the trajectory
does not return to the domain $\mathbb{Z}_{p}$ in finite time.
Thus, for $\alpha<1$, the random walk is nonrecurrent.

For $\alpha\ge1$ we can represent the process of random walk on
$\mathbb{Q}_{p}$ as an infinite sequence of events consisting in
successive sojourns of the trajectory in the sets
$\mathbb{Z}_{p}$, $\mathbb{Q}_{p}\setminus\mathbb{Z}_{p}$,
$\mathbb{Z}_{p}$, $\mathbb{Q}_{p}\setminus\mathbb{Z}_{p}$,
$\ldots$ with random independent sojourn times
$\zeta_{1},\:\eta_{1},\:\zeta_{2},\:\eta_{2},\:\ldots$, where
$\zeta_{1},\:\zeta_{2},\ldots$, just as
$\eta_{1},\:\eta_{2},\ldots$, have the same distribution. In this
case, the distribution function $\zeta_{i}$ and its density are
given by

\begin{equation}
G\left(t\right)=\Theta\left(t\right)\left(1-\exp\left(-B_{\alpha}t\right)\right),\;g\left(t\right)=\dfrac{dG\left(t\right)}{dt}=B_{\alpha}\exp\left(-B_{\alpha}t\right),\label{G_t}
\end{equation}
respectively. In terms of Laplace transforms, the functions
(\ref{G_t}) are rewritten as
\begin{equation}
\hat{G}\left(s\right)=\dfrac{1}{s}\hat{g}\left(s\right),\:\hat{g}\left(s\right)=\dfrac{B_{\alpha}}{s+B_{\alpha}}.\label{G_s}
\end{equation}
The distribution function $\eta_{i}$ is
\[
H\left(t\right)=\intop_{0}^{t}h\left(t\right)dt
\]
and is explicitly defined by the Laplace transform
$\hat{H}\left(s\right)=\dfrac{1}{s}\hat{h}\left(s\right)$, where
$\hat{h}\left(s\right)$ is the Laplace transform of the
distribution function density $h\left(t\right)$, which, by the
convolution relations

\[
f\left(t\right)=\intop_{0}^{t}g\left(t-\theta\right)h\left(\theta\right)d\theta
\]
is expressed as
\begin{equation}
\hat{h}\left(s\right)=\dfrac{\hat{f}\left(s\right)}{\hat{g}\left(s\right)}=B_{\alpha}^{-1}\left(s+B_{\alpha}-\frac{1}{\hat{J}\left(s\right)}\right).\label{h_s}
\end{equation}
The sum of random variables $\alpha_{n}=\sum_{i=1}^{n}\zeta_{i}$
is the total sojourn time in the set $\mathbb{Z}_{p}$ provided
that the trajectory hits $\mathbb{Z}_{p}$ precisely $n$ times. In
this case, the distribution function $G_{n}\left(t\right)$ of the
sum $\alpha_{n}$ is determined by the Laplace transform

\begin{equation}
\hat{G}_{n}\left(s\right)=\dfrac{1}{s}\left(\dfrac{B_{\alpha}}{s+B_{\alpha}}\right)^{n},\label{G_n_s}
\end{equation}
whose inverse is

\begin{equation}
G_{n}\left(t\right)=1-\sum_{m=0}^{n-1}\dfrac{B_{\alpha}^{m}}{m!}t^{m}\exp\left(-B_{\alpha}t\right).\label{G_n_t}
\end{equation}

In the case of a recurrent random walk ($\alpha\geq1$), the sum of
random variables
$t_{\mathbb{Q}_{p}\setminus\mathbb{Z}_{p},n}\equiv\sum_{i=1}^{n}\eta_{i}$
is the total sojourn time in the set
$\mathbb{Q}_{p}\setminus\mathbb{Z}_{p}$ provided that the
trajectory hits $\mathbb{Q}_{p}\setminus\mathbb{Z}_{p}$ precisely
$n$ times. In this case, the distribution function
$H_{n}\left(t\right)$ of the sum
$t_{\mathbb{Q}_{p}\setminus\mathbb{Z}_{p},n}$ is determined by
the Laplace transform
\begin{equation}
\hat{H}_{n}\left(s\right)=\dfrac{1}{s}\hat{h}_{n}\left(s\right),\label{H_n_s}
\end{equation}
where the Laplace transform $\hat{h}_{n}\left(s\right)$ of the
density $h_{n}\left(t\right)$ is given by

\begin{equation}
\hat{h}_{n}\left(s\right)=B_{\alpha}^{-n}\left(s+B_{\alpha}-\frac{1}{\hat{J}\left(s\right)}\right)^{n}.\label{h_n_s}
\end{equation}

For $\alpha<1$ we can also represent the process of a random walk
on $\mathbb{Q}_{p}$ as a successive alternate sojourn of the
trajectory in the sets
$\mathbb{Z}_{p}$$,\:\mathbb{Q}_{p}\setminus\mathbb{Z}_{p},\:\mathbb{Z}_{p},\:\mathbb{Q}_{p}\setminus\mathbb{Z}_{p},\ldots$
with independent random sojourn times, defined by the sequence
$\zeta_{1},\:\eta_{1},\:\zeta_{2},\:\eta_{2},\:\ldots$. In this
case, we should assume that $\eta_{1},\:\eta_{2},\ldots$ take
values on the extended nonnegative number axis
$\mathbb{\bar{R}}_{+}=\mathbb{R}_{+}\cup\left\{ +\infty\right\} $.
Note that, for any specific implementation of this process, the
sequence $\zeta_{1},\:\eta_{1},\:\zeta_{2},\:\eta_{2},\:\ldots$ is
truncated at some number $i=k$ such that the value of
$\eta_{k}$ on the extended number axis is equal
to $+\infty$. In this case, the distribution function
$G_{n}\left(t\right)$ is also determined by the formula
(\ref{G_n_t}), and the distribution function $H\left(t\right)$ is
\[
H\left(t\right)=\left\{ \begin{array}{c}
\intop_{0}^{t}h\left(t\right)dt,\:t<+\infty\\
1,\:t=+\infty
\end{array}\right.,
\]
where $h\left(t\right)$ is defined by the Laplace transform
(\ref{h_s}). Note that in this case the probability of exit of the
trajectory from $\mathbb{Q}_{p}\setminus\mathbb{Z}_{p}$ is
different from $1$ and equals

\[
\intop_{0}^{+\infty}h\left(t\right)dt=\hat{h}\left(0\right)=\hat{f}\left(0\right)=C_{\alpha}<1\:\mathrm{with}\:\alpha<1.
\]

When $\alpha<1$, we can also introduce a random variable
$\beta_{n}$ -- the total sojourn time in the set
$\mathbb{Q}_{p}\setminus\mathbb{Z}_{p}$ provided that the
trajectory hits $\mathbb{Q}_{p}\setminus\mathbb{Z}_{p}$ precisely
$n$ times, taking values on the extended nonnegative number axis
$\mathbb{\bar{R}}_{+}$. Note that in this case $\beta_{n}$ cannot
be represented as a sum of independent random variables
$\eta_{i}$. The distribution function $H_{n}\left(t\right)$ of the
random variable $\beta_{n}$ is equal to
\[
H_{n}\left(t\right)=\left\{ \begin{array}{c}
\intop_{0}^{t}h_{n}\left(t\right)dt,\:t<+\infty\\
1,\:t=+\infty
\end{array}\right.,
\]
where $h_{n}\left(t\right)$ is defined by its Laplace transform
(\ref{h_n_s}). Note that
$\intop_{0}^{\infty}h_{n}\left(t\right)dt=C_{\alpha}^{n}<1.$

\section{Distribution of the sojourn time in the domain $\mathbb{Z}_{p}$ }

Denote the integral distribution functions of the  processes
(\ref{th_Z_p}) and (\ref{th_Q_p-Z_p}) by
$\Phi_{\mathbb{Z}_{p}}\left(\theta,t\right)$ and
$\Phi_{\mathbb{Q}_{p}\setminus\mathbb{Z}_{p}}\left(\theta,t\right)$,
respectively, and their densities, by
$\phi_{\mathbb{Z}_{p}}\left(\theta,t\right)$ and
$\phi_{\mathbb{Q}_{p}\setminus\mathbb{Z}_{p}}\left(\theta,t\right)$,
respectively. Since
\begin{equation}
\theta_{\mathbb{Z}_{p}}\left(t\right)+\theta_{\mathbb{Q}_{p}\setminus\mathbb{Z}_{p}}\left(t\right)=t,\label{th+th}
\end{equation}
it follows that
$\phi_{\mathbb{Z}_{p}}\left(\theta,t\right)=\phi_{\mathbb{Q}_{p}\setminus\mathbb{Z}_{p}}\left(t-\theta,t\right)$
and the functions $\Phi_{\mathbb{Z}_{p}}\left(\theta,t\right)$ and
$\Phi_{\mathbb{Q}_{p}\setminus\mathbb{Z}_{p}}\left(\theta,t\right)$
are related by
\begin{equation}
\Phi_{\mathbb{Z}_{p}}\left(\theta,t\right)=1-\Phi_{\mathbb{Q}_{p}\setminus\mathbb{Z}_{p}}\left(t-\theta,t\right).\label{Phi_z}
\end{equation}

\textbf{Theorem 3}. The distribution functions
$\Phi_{\mathbb{Q}_{p}\setminus\mathbb{Z}_{p}}\left(\theta,t\right)$
and $\Phi_{\mathbb{Z}_{p}}\left(\theta,t\right)$ are given by

\begin{equation}
\Phi_{\mathbb{Q}_{p}\setminus\mathbb{Z}_{p}}\left(\theta,t\right)=\sum_{n=0}^{\infty}H_{n}\left(\theta\right)\left(G_{n}\left(t-\theta\right)-G_{n+1}\left(t-\theta\right)\right),\label{Phi_B}
\end{equation}
\begin{equation}
\Phi_{\mathbb{Z}_{p}}\left(\theta,t\right)=1-\sum_{n=0}^{\infty}H_{n}\left(t-\theta\right)\left(G_{n}\left(\theta\right)-G_{n+1}\left(\theta\right)\right).\label{Phi_A}
\end{equation}
where $G_{n}\left(t\right)$ ($n>0$) is defined by formula
(\ref{G_n_t}), $H_{n}\left(t\right)$ ($n>0$) is defined by the
Laplace transforms (\ref{H_n_s}) -- (\ref{h_n_s}), and for $n=0$
we set $G_{0}\left(t\right)\equiv1$ and
$H_{0}\left(t\right)=\Theta\left(t\right)$.

\textbf{Proof.} This theorem was proved in the general case in
\cite{Takacs_1,Takacs_2,Takacs_3} by two methods. In view of the
importance of this result, we reproduce one of the proofs in full,
with its adaptation to our case of a $p$-adic random walk. The
proof given here covers the cases of recurrent
$\left(\alpha\geq1\right)$ and nonrecurrent
$\left(\alpha<1\right)$ random walks.

Let us fix $t$. We also fix a number $\theta$, $0\leq\theta<t$,
and consider the event
$\theta_{\mathbb{Q}_{p}\setminus\mathbb{Z}_{p}}\left(t\right)\leq\theta$.
Next, consider the equation
\begin{equation}
\theta_{\mathbb{Z}_{p}}\left(u\right)=t-\theta\label{u}
\end{equation}
for $u$. Denote by $\tau=\tau\left(t-\theta\right)$ the random
variable that is the smallest of all the solutions
$u\in\left[0,+\infty\right)$ satisfying the equation (\ref{u}).
For fixed $t$ and $\theta$, $\tau$ is a random variable. Let us
determine the existence condition of $\tau$.

For a recurrent random walk ($\alpha\geq1$), the solution $\tau$
exists always. One can verify this by the following arguments.
Consider a specific implementation of a random process and suppose
that the trajectory visits $\mathbb{Z}_{p}$ precisely $n+1$ times
during the time from $0$ to $u$. Denote
$t_{\mathbb{Z}_{p},n}=\zeta_{1}+\zeta_{2}+\cdots\zeta_{n}$ and
$t_{\mathbb{Q}_{p}\setminus\mathbb{Z}_{p},n}=\eta_{1}+\eta_{2}+\cdots\eta_{n}$.
In this case, $t_{\mathbb{Z}_{p},n}+\Delta t=t-\theta$, where
$\Delta t<\zeta_{n+1}$ (if the system is in $\mathbb{Z}_{p}$ at
time $u$) or $\Delta t=\zeta_{n+1}$ (if the system is in
$\mathbb{Q}_{p}\setminus\mathbb{Z}_{p}$ at time $u$). In the first
case, the solution of equation (\ref{u}) is unique, and $\tau=u$.
In the second case, equation (\ref{u}) has an infinite number of
solutions that differ by the time of the last sojourn in
$\mathbb{Q}_{p}\setminus\mathbb{Z}_{p}$, and then $\tau$ is the
smallest of all solutions $u$ satisfying equation (\ref{u}), where
$\theta_{\mathbb{Z}_{p}}\left(\tau\right)=t-\theta$,
$\xi\left(\tau\right)\in\mathbb{Z}_{p}$. Thus, $\tau$ is a current
instant of time if the trajectory visits $\mathbb{Z}_{p}$ the
($n+1$)th time (in the first case) or the time of the first
transition
$\mathbb{Z}_{p}\rightarrow\mathbb{Q}_{p}\setminus\mathbb{Z}_{p}$
after the system visited $\mathbb{Z}_{p}$ the ($n+1$)th time (in
the second case). It is obvious that
$\theta_{\mathbb{Q}_{p}\setminus\mathbb{Z}_{p}}\left(\tau\right)=t_{\mathbb{Q}_{p}\setminus\mathbb{Z}_{p},n}$
and
$t_{\mathbb{Z}_{p},n}<\theta_{\mathbb{Z}_{p}}\left(\tau\right)\leq
t_{\mathbb{Z}_{p},n+1}$ .

For a nonrecurrent random walk ($\alpha<1$), an event is possible
in which the trajectory hits the set
$\mathbb{Q}_{p}\setminus\mathbb{Z}_{p}$ $n$ times and then does
not leave it in finite time. In this case, for sufficiently small
$\theta$, equation (\ref{u}) may not have a solution for $u$ in
the domain $u\in\left[0,+\infty\right)$. Then, for any
$u\in\left[0,+\infty\right)$ we have
$\theta_{\mathbb{Z}_{p}}\left(u\right)<t-\theta$. Then, by
(\ref{th+th}), it follows that
$\theta_{\mathbb{Z}_{p}}\left(u\right)<\theta_{\mathbb{Z}_{p}}\left(t\right)+\theta_{\mathbb{Q}_{p}\setminus\mathbb{Z}_{p}}\left(t\right)-\theta$;
for $u=t$, this implies
$\theta_{\mathbb{Q}_{p}\setminus\mathbb{Z}_{p}}\left(t\right)>\theta$.
If
$\theta_{\mathbb{Q}_{p}\setminus\mathbb{Z}_{p}}\left(t\right)\leq\theta$,
then $\theta_{\mathbb{Z}_{p}}\left(t\right)>t-\theta$, and the
random variable $\tau$ exists.

Thus, the variable $\tau$ always exists in the probability
subspace defined by the condition
$\theta_{\mathbb{Q}_{p}\setminus\mathbb{Z}_{p}}\left(t\right)\leq\theta$.

Next, taking into account the equalities
$\theta_{\mathbb{Z}_{p}}\left(\tau\right)=t-\theta$ and
$\theta_{\mathbb{Z}_{p}}\left(t\right)+\theta_{\mathbb{Q}_{p}\setminus\mathbb{Z}_{p}}\left(t\right)=t$
and the fact that $\theta_{\mathbb{Z}_{p}}\left(t\right)$ and
$\theta_{\mathbb{Q}_{p}\setminus\mathbb{Z}_{p}}\left(t\right)$ are
nondecreasing functions of $t$, we have a chain of equivalent
events:

\[
\theta_{\mathbb{Q}_{p}\setminus\mathbb{Z}_{p}}\left(t\right)\leq\theta\:\Leftrightarrow\:\theta_{\mathbb{Z}_{p}}\left(\tau\right)\leq\theta_{\mathbb{Z}_{p}}\left(t\right)\:\Leftrightarrow\:\tau\leq t\:\Leftrightarrow\:\theta_{\mathbb{Z}_{p}}\left(\tau\right)+\theta_{\mathbb{Q}_{p}\setminus\mathbb{Z}_{p}}\left(\tau\right)\leq t\:\Leftrightarrow\:\theta_{\mathbb{Q}_{p}\setminus\mathbb{Z}_{p}}\left(\tau\right)\leq\theta
\]
Therefore, the probabilities of the events
$\theta_{\mathbb{Q}_{p}\setminus\mathbb{Z}_{p}}\left(t\right)\leq\theta$
and
$\theta_{\mathbb{Q}_{p}\setminus\mathbb{Z}_{p}}\left(\tau\right)\leq\theta$
are identical:
\[
\mathrm{P}\left[\theta_{\mathbb{Q}_{p}\setminus\mathbb{Z}_{p}}\left(t\right)\leq\theta\right]=\mathrm{P}\left[\theta_{\mathbb{Q}_{p}\setminus\mathbb{Z}_{p}}\left(\tau\right)\leq\theta\right].
\]
Notice that, for a fixed $n$, the equality
$\theta_{B}\left(\tau\right)=t_{\mathbb{Q}_{p}\setminus\mathbb{Z}_{p},n}$
holds provided that $t_{\mathbb{Z}_{p},n}<t-\theta\leq
t_{\mathbb{Z}_{p},n+1}$. Thus, the event
$\theta_{\mathbb{Q}_{p}\setminus\mathbb{Z}_{p}}\left(\tau\right)\leq\theta$
is decomposed into a sum of $n=0,1,2,\ldots$ independent events
$\left(t_{\mathbb{Q}_{p}\setminus\mathbb{Z}_{p},n}\leq\theta\right)\cup\left(t_{\mathbb{Z}_{p},n}<t-\theta\leq
t_{\mathbb{Z}_{p},n+1}\right)$. Then

\[
\mathrm{P}\left[\theta_{\mathbb{Q}_{p}\setminus\mathbb{Z}_{p}}\left(t\right)\leq\theta\right]=\mathrm{P}\left[\theta_{\mathbb{Q}_{p}\setminus\mathbb{Z}_{p}}\left(\tau\right)\leq\theta\right]=\sum_{n=0}^{\infty}\mathrm{P}\left[\left(t_{\mathbb{Q}_{p}\setminus\mathbb{Z}_{p},n}\leq\theta\right)\cup\left(t_{Z_{p},n}<t-\theta\leq
t_{\mathbb{Z}_{p},n+1}\right)\right];
\]
this implies (\ref{Phi_B}). Then (\ref{Phi_z}) implies
(\ref{Phi_A}), which proves Theorem 3.

\section{The mean sojourn time in $\mathbb{Z}_{p}$ and the asymptotics as $t\rightarrow\infty$
of the moments of the distribution of the sojourn time in
$\mathbb{Z}_{p}$.
}

\textbf{Theorem 4. } The mean sojourn time of a trajectory in
$\mathbb{Z}_{p}$ is defined by
\begin{equation}
\left\langle \theta\right\rangle \left(t\right)=\intop_{0}^{t}J\left(t'\right)dt'\label{theta_mid}
\end{equation}
\textbf{Proof.}

Applying formula (\ref{Phi_A}) and integrating by parts with
regard to $\Phi_{\mathbb{Z}_{p}}\left(0,t\right)=0$ and
$\Phi_{\mathbb{Z}_{p}}\left(t,t\right)=1$, we obtain

\begin{equation}
\left\langle \theta\right\rangle \left(t\right)=\intop_{0}^{t}\theta d\Phi_{\mathbb{Z}_{p}}\left(\theta,t\right)=\sum_{n=0}^{\infty}\intop_{0}^{t}H_{n}\left(t-\theta\right)\left(G_{n}\left(\theta\right)-G_{n+1}\left(\theta\right)\right).\label{M_th}
\end{equation}
Take the Laplace transform of (\ref{M_th}). Then, taking into
account (\ref{G_n_s}) and (\ref{h_n_s}), we obtain

\[
\hat{\left\langle \theta\right\rangle }\left(s\right)=\sum_{n=0}^{\infty}\hat{H}_{n}\left(s\right)\left(\hat{G}_{n}\left(s\right)-\hat{G}_{n+1}\left(s\right)\right)=
\]
\[
=\dfrac{1}{s}\sum_{n=0}^{\infty}\left(\left(s+B_{\alpha}\right)-\frac{1}{J(s)}\right)^{n}\left(\dfrac{1}{s+B_{\alpha}}\right)^{n}\dfrac{1}{s+B_{\alpha}}=
\]
\[
=\dfrac{1}{s}\dfrac{1}{s+B_{\alpha}-\left(\left(s+B_{\alpha}\right)-\frac{1}{J(s)}\right)}=\dfrac{1}{s}\hat{J}(s).
\]
Passing to the Laplace originals, we obtain (\ref{theta_mid}),
which proves the theorem.

Next, we establish the asymptotic behavior of the moments
$\left\langle \theta^{\beta}\right\rangle \left(t\right)$ as
$t\rightarrow\infty$ for $\beta>0$. To this end, we need several
theorems.

\textbf{Theorem 5.} Let
\[
\Psi\left(\theta,t\right)\equiv\sum_{n=0}^{\infty}\tilde{H}_{n}\left(\theta\right)\left(G_{n}\left(t-\theta\right)-G_{n+1}\left(t-\theta\right)\right),
\]
where $G_{n}\left(t\right)$ is the $n$ fold iterated convolution
of some distribution function $G\left(t\right)$ that has a finite
first moment $\tilde{\alpha}$ and a finite second moment
$\tilde{\beta}$, and $\tilde{H}_{n}\left(\theta\right)$ is the $n$
fold iterated convolution of some distribution function
$\tilde{H}\left(\theta\right)$ such that
\begin{equation}
B=\lim_{t\rightarrow\infty}\left(1-\tilde{H}\left(t\right)\right)t^{\gamma}\label{B}
\end{equation}
is a finite quantity for some $0<\gamma<1$,
$G_{0}\left(t\right)\equiv1$, and
$H_{0}\left(t\right)=\Theta\left(t\right)$. Then, for any $x>0$,

\[
\lim_{t\rightarrow\infty}\Psi\left(t-xt^{\gamma},t\right)=F_{\gamma}\left(\left(\dfrac{\tilde{\alpha}}{Bx}\right)^{\tfrac{1}{\gamma}}\right),
\]
where $F_{\gamma}\left(y\right)$ is the integral function of
stable distribution having the characteristic function
\begin{equation}
\tilde{f}_{\gamma}\left(k\right)=\exp\left(-\left|k\right|^{\gamma}\left(\cos\left(\dfrac{\pi\gamma}{2}\right)-i\sin\left(\dfrac{\pi\gamma}{2}\right)\mathrm{sign}\left(k\right)\right)\Gamma\left(1-\gamma\right)\right).\label{har_stable}
\end{equation}
The proof of Theorem 5 is given in \cite{Takacs_4}.

It follows from Theorem 5 that for the function
\[
\tilde{\Psi}\left(\theta,t\right)=1-\Psi\left(t-\theta,t\right)
\]
we have

\begin{equation}
\lim_{t\rightarrow\infty}\tilde{\Psi}\left(xt^{\gamma},t\right)=1-F_{\gamma}\left(\left(\dfrac{1}{B_{\alpha}Bx}\right)^{\tfrac{1}{\gamma}}\right).\label{Psi_til}
\end{equation}
Since the function
$\tilde{\Psi}\left(xt^{\gamma},t\right)=1-\Psi\left(t-xt^{\gamma},t\right)$
is a probability measure by its meaning, it is bounded for $x>0$.
Therefore, the limit (\ref{Psi_til}) converges uniformly in $x$,
and we have

\begin{equation}
\dfrac{d}{dx}\lim_{t\rightarrow\infty}\tilde{\Psi}\left(xt^{\gamma},t\right)=\lim_{t\rightarrow\infty}t^{\gamma}\tilde{\psi}\left(xt^{\gamma},t\right)=f_{\gamma}\left(\left(\dfrac{1}{B_{\alpha}Bx}\right)^{\tfrac{1}{\gamma}}\right)\dfrac{1}{\gamma}\left(\dfrac{1}{B_{\alpha}Bx}\right)^{\tfrac{1}{\gamma}}x^{\tfrac{-\gamma-1}{\gamma}},\label{Psi_til_diff}
\end{equation}
where $f_{\gamma}\left(t\right)$ is the density of the function of
stable distribution with the characteristic function
(\ref{har_stable}) and
$\tilde{\psi}\left(\theta,t\right)=\dfrac{d}{d\theta}\tilde{\Psi}\left(\theta,t\right)$.

The distribution (\ref{G_t}) has finite first two moments, namely,
$\left\langle t_{\mathbb{Z}_{p},1}\right\rangle
=\dfrac{1}{B_{\alpha}}\equiv\tilde{\alpha}$ and $\left\langle
t_{\mathbb{Z}_{p},1}^{2}\right\rangle
=\dfrac{1}{B_{\alpha}^{2}}\equiv\tilde{\beta}$. In addition, for
$\alpha\geq1$, Theorem 2 implies that

\begin{equation}
C_{1}\le\left(1-H(t)\right)t^{\gamma}\le
C_{2},\:\alpha>1,\:\gamma=\dfrac{\alpha-1}{\alpha},\label{1-H}
\end{equation}
where $C_{i}$, $i=1,2$, are some functions of $\alpha$ and $p$.
This gives us right to apply formula (\ref{Psi_til_diff}) for the
function $\phi_{\mathbb{Z}_{p}}\left(xt^{\gamma},t\right)$.
Namely, taking into account (\ref{1-H}) and (\ref{Psi_til}), we
obtain the following equation  also for $\alpha>1$ also for
$\alpha>1$ for $x>0$:

\begin{equation}
g_{\min}\left(x\right)\leq\lim_{t\rightarrow\infty}t^{\gamma}\phi_{\mathbb{Z}_{p}}\left(xt^{\gamma},t\right)\leq g_{\max}\left(x\right),\label{noneq_for_f}
\end{equation}
where
\[
g_{\min}\left(x\right)=\min_{i=1,2}\left\{ f_{\gamma}\left(\left(\dfrac{1}{B_{\alpha}C_{i}x}\right)^{\tfrac{1}{\gamma}}\right)\left(\dfrac{1}{B_{\alpha}C_{i}x}\right)^{\tfrac{1}{\gamma}}\dfrac{1}{\gamma}x^{\tfrac{-\gamma-1}{\gamma}}\right\} ,
\]
\[
g_{\max}\left(x\right)=\max_{i=1,2}\left\{ f_{\gamma}\left(\left(\dfrac{1}{B_{\alpha}C_{i}x}\right)^{\tfrac{1}{\gamma}}\right)\left(\dfrac{1}{B_{\alpha}C_{i}x}\right)^{\tfrac{1}{\gamma}}\dfrac{1}{\gamma}x^{\tfrac{-\gamma-1}{\gamma}}\right\} .
\]
Unfortunately, for $\alpha\leq1$ we cannot apply Theorem 5 to
obtain an estimate similar to (\ref{noneq_for_f}), because

\[
C_{1}\left(\log t\right)^{-1}\left(1+o(1)\right)\le1-H(t)\le C_{2}\left(\log t\right)^{-1}\left(1+o(1)\right),\:\alpha=1,
\]
\[
C_{1}t^{-\gamma}\left(1+o(1)\right)\le C_{\alpha}-H(t)\le C_{2}t^{-\gamma}\left(1+o(1)\right),\:\gamma=\dfrac{1-\alpha}{\alpha}\:\alpha<1,
\]
and condition (\ref{B}) cannot be satisfied.

\textbf{Lemma 1. } For any $\varepsilon>0$, the density of the
function of stable distribution (\ref{har_stable}) admits a
representation in the form of a series uniformly convergent in
$t\geq\varepsilon$,

\begin{equation}
f_{\gamma}\left(t\right)=\dfrac{1}{\pi}\sum_{n=1}^{\infty}\left(-1\right)^{n+1}\dfrac{\Gamma\left(1-\gamma\right)^{n}\sin\left(n\pi\gamma\right)}{n!}\Gamma\left(n\gamma+1\right)t^{-n\gamma-1}.\label{Lemma_f_gamma}
\end{equation}

\textbf{Proof. } Denoting where $a=\Gamma\left(1-\gamma\right)$,
we write

\[
f_{\gamma}\left(t\right)=\dfrac{1}{2\pi}\intop_{-\infty}^{+\infty}\exp\left(-ikt\right)\exp\left(-\left|k\right|^{\gamma}\left(\cos\left(\dfrac{\pi\gamma}{2}\right)-i\sin\left(\dfrac{\pi\gamma}{2}\right)\mathrm{sign}\left(k\right)\right)a\right)ds=
\]
\[
=\dfrac{1}{\pi}\mathrm{Re}\intop_{0}^{+\infty}\exp\left(-ikt\right)\exp\left(-a\exp\left(-i\dfrac{\pi\gamma}{2}\right)k^{\gamma}\right)dk
\]
Applying the Cauchy theorem, we can represent this integral in the form
\[
\ointop_{C}\exp\left(-izt\right)\exp\left(-\exp\left(-i\dfrac{\pi\gamma}{2}\right)z^{\gamma}\right)dz
\]
where the contour $C$ is bounded by the quarters of two circles
with $\mathrm{Re}z>0$, $\mathrm{Im}z$\textless 0 and radii $r$ and
$R$, respectively, and two segments $\left[r,R\right]$ on the real
and imaginary axes. Passing to the limit as $r\rightarrow0$ and
$R\rightarrow\infty$ and applying the Jordan lemma, we obtain
\[
\intop_{0}^{+\infty}\exp\left(-ikt\right)\exp\left(-a\exp\left(-i\dfrac{\pi\gamma}{2}\right)k^{\gamma}\right)dk=-i\intop_{0}^{+\infty}\exp\left(-kt\right)\exp\left(-a\exp\left(-i\pi\gamma\right)k^{\gamma}\right)dk;
\]
hence,
\begin{equation}
f_{\gamma}\left(t\right)=-\dfrac{1}{\pi}\mathrm{Re}i\intop_{0}^{+\infty}\exp\left(-kt\right)\exp\left(-a\exp\left(-i\pi\gamma\right)k^{\gamma}\right)dk.\label{f_gamma}
\end{equation}
Expanding the exponent in a series  and changing the variable
$kt=x$, we obtain
\[
f_{\gamma}\left(t\right)=-\dfrac{1}{\pi t}\mathrm{Re}i\sum_{n=0}^{\infty}\left(-1\right)^{n}\dfrac{a^{n}\exp\left(-ni\pi\gamma\right)}{n!}t^{-n\gamma}\intop_{0}^{+\infty}\exp\left(-x\right)x^{n\gamma}dx=
\]
\[
=\dfrac{1}{\pi}\sum_{n=1}^{\infty}\left(-1\right)^{n+1}\dfrac{a^{n}\sin\left(n\pi\gamma\right)}{n!}\Gamma\left(n\gamma+1\right)t^{-n\gamma-1},
\]
which implies the assertion of the lemma.

\textbf{Theorem 6}.\textbf{ }For $\alpha\text{>}1$ and
$t\rightarrow\infty,$
\begin{equation}
D_{\min}t^{\tfrac{\alpha-1}{\alpha}\beta}\left(1+o\left(1\right)\right)\leq\left\langle
\theta^{\beta}\right\rangle \left(t\right)\leq
D_{\max}t^{\tfrac{\alpha-1}{\alpha}\beta}\left(1+o\left(1\right)\right),\:\beta<\dfrac{\alpha}{\alpha-1},\label{T7_1}
\end{equation}
\begin{equation}
\tilde{D}_{\min}t^{\beta-\tfrac{1}{\alpha-1}}\left(1+o\left(1\right)\right)\leq\left\langle
\theta^{\beta}\right\rangle
\left(t\right)\leq\tilde{D}_{\max}t^{\beta-\tfrac{1}{\alpha-1}}\left(1+o\left(1\right)\right),\:\beta\geq\dfrac{\alpha}{\alpha-1},\label{T7_2}
\end{equation}
where $D_{\min},\:D_{\max},\:\tilde{D}_{\min},\:\tilde{D}_{\max}$
are some functions of $\alpha$ and $\beta$.

\textbf{Proof.}

Denote $\dfrac{\alpha-1}{\alpha}\equiv\gamma$. Consider
\[
\left\langle \theta^{\beta}\right\rangle
\left(t\right)=\intop_{0}^{t}\theta^{\beta}\dfrac{d\Phi_{\mathbb{Z}_{p}}\left(\theta,t\right)}{d\theta}d\theta=\intop_{0}^{t}\theta^{\beta}\phi_{\mathbb{Z}_{p}}\left(\theta,t\right)d\theta.
\]
Let us change the variable $\theta=xt^{\gamma}$:
\[
\left\langle \theta^{\beta}\right\rangle \left(t\right)=t^{\left(\beta+1\right)\gamma}\intop_{0}^{t^{1-\gamma}}x^{\beta}\phi_{\mathbb{Z}_{p}}\left(xt^{\gamma},t\right)dx.
\]
By (\ref{noneq_for_f}),
\[
M_{\min}\left(t\right)\left(1+o\left(1\right)\right)\leq\left\langle \theta^{\beta}\right\rangle \left(t\right)\leq M_{\max}\left(t\right)\left(1+o\left(1\right)\right),
\]
\begin{equation}
M_{\min}\left(t\right)=\dfrac{1}{\gamma}\left(\dfrac{1}{B_{\alpha}C_{i}}\right)^{\tfrac{1}{\gamma}}t^{\beta\gamma}\intop_{0}^{t^{1-\gamma}}x^{\beta-1}x^{-\tfrac{1}{\gamma}}\min_{i=1,2}\left\{ \left(\dfrac{1}{B_{\alpha}C_{i}}\right)^{\tfrac{1}{\gamma}}f_{\gamma}\left(\left(\dfrac{1}{B_{\alpha}C_{i}x}\right)^{\tfrac{1}{\gamma}}\right)\right\} dx,\label{M_min}
\end{equation}
\begin{equation}
M_{\max}\left(t\right)=\dfrac{1}{\gamma}t^{\beta\gamma}\intop_{0}^{t^{1-\gamma}}x^{\beta-1}x^{-\tfrac{1}{\gamma}}\max_{i=1,2}\left\{ \left(\dfrac{1}{B_{\alpha}C_{i}}\right)^{\tfrac{1}{\gamma}}f_{\gamma}\left(\left(\dfrac{1}{B_{\alpha}C_{i}x}\right)^{\tfrac{1}{\gamma}}\right)\right\} dx.\label{M_max}
\end{equation}
By Lemma 1 we have
$f_{\gamma}\left(\left(\dfrac{1}{B_{\alpha}C_{i}x}\right)^{\tfrac{1}{\gamma}}\right)=A_{i}x^{\tfrac{\gamma+1}{\gamma}}\left(1+O\left(x\right)\right)$
as $x\rightarrow0$; therefore, the integrands in (\ref{M_min}) and
(\ref{M_max}) behave as $x^{\beta}$ as $x\rightarrow0$. Since
$f_{\gamma}\left(0\right)$ is finite, as $x\rightarrow\infty$, the
integrand behaves as
$x^{\beta-1-\tfrac{1}{\gamma}}=x^{\beta-1-\tfrac{\alpha}{\alpha-1}}$.
Therefore, for $i=1,2$, the integral
\[
\intop_{0}^{\infty}x^{\beta-1}x^{-\tfrac{1}{\gamma}}f_{\gamma}\left(\left(\dfrac{1}{B_{\alpha}C_{i}x}\right)^{\tfrac{1}{\gamma}}\right)dx
\]
converges for $\beta<\dfrac{\alpha}{\alpha-1}$, and in this case
we have, as $t\rightarrow\infty$,

\[
M_{\min}\left(t\right)=D_{\min}t^{\beta\gamma}=D_{\min}t^{\tfrac{\alpha-1}{\alpha}\beta},
\]
\[
M_{\max}\left(t\right)=D_{\max}t^{\beta\gamma}=D_{\max}t^{\tfrac{\alpha-1}{\alpha}\beta},
\]
for some $D_{\min}$ and $D_{\max}$. To find the asymptotics as
$t\rightarrow\infty$ in the case of
$\beta\geq\dfrac{\alpha}{\alpha-1}$, we also should take into
account the contribution of the upper limit of the integrals in
(\ref{M_min}) and (\ref{M_max})
\[
t^{\beta\gamma}\intop_{0}^{t^{1-\gamma}}x^{\beta-1}x^{-\tfrac{1}{\gamma}}f_{\gamma}\left(\left(\dfrac{1}{B_{\alpha}C_{i}x}\right)^{\tfrac{1}{\gamma}}\right)dx=
\]

\[
=t^{\beta\gamma}\intop_{\epsilon}^{t^{1-\gamma}}x^{\beta-1-\tfrac{1}{\gamma}-1}f_{\gamma}\left(0\right)dx\left(1+o\left(1\right)\right)=
\]
\[
=\dfrac{1}{\beta-\tfrac{1}{\gamma}}t^{\beta\gamma}t^{\left(\beta-\tfrac{1}{\gamma}\right)\left(1-\gamma\right)}f_{\gamma}\left(0\right)dx\left(1+o\left(1\right)\right)=
\]
\[
=\dfrac{1}{\beta-\tfrac{\alpha}{\alpha-1}}t^{\beta-\tfrac{1}{\alpha-1}}f_{\gamma}\left(0\right)dx\left(1+o\left(1\right)\right),
\]
which proves the theorem.

\section{Conclusions and possible applications}

The main result of this work is the asymptotic behavior as
$t\rightarrow\infty$ of the moments of the distribution of a
random variable -- the sojourn time of the trajectory in
$\mathbb{Z}_{p}$ for a $p$-adic random walk for $\alpha>1$. This
result is formulated in Theorem 6. Here we discuss one of possible
applications of this result.

As already mentioned in the Introduction, one of the most
important applications of the process of $p$-adic random walk is
modeling the conformational dynamics of protein molecules in the
native state. Of particular interest are processes in protein that
can occur only if the protein is in a certain domain of the
conformational space of states. It is known that, in many
experiments, the fluctuation dynamics of a protein molecule can be
investigated by following the random variation of the optical
absorption frequency of a chromophore marker located in the active
center of the protein. This frequency of optical absorption of the
marker is sensitive to the spatial arrangement of the surrounding
atoms. In experiments, one often applies the method of ``frequency
coloring'' of a small part of markers in a sample by inducing an
irreversible photochemical transition in it by pulsed
monochromatic pumping. As a result, in subsequent б­пвЁЁ of the
absorption spectrum of a sample, this part of markers turns out to
be ``dark,'' and a spectral hole with the characteristic width on
the order of the absorption bandwidth of the marker is formed in
the absorption spectrum. This techniques is called the method of
spectral hole burnout. When the absorption frequency of the
markers located in individual protein molecules of a sample
randomly changes due to the rearrangements of their atomic
surrounding, the spectral hole broadens. This phenomenon is called
spectral diffusion. By the change in the width of the spectral
hole with time, one can judge of the fluctuation dynamics of the
structure under study.

Spectral diffusion was studied at low temperatures for various
globular proteins embedded in different organic matrices (see
\cite{LW,PFVBB,Friedrich1,Friedrich2} and references cited
therein). The properties of spectral diffusion in proteins were
determined by the variation of the characteristic width of the
spectral hole $\sigma_{\nu}$ as a function of time $t$
since the burnout time (the ``observation'' time). In addition, we
have investigated the dependence of the characteristic width of
the spectral hole $\sigma_{\nu}$ on the interval between
the preparation time of the sample and the burnout time of the
spectral hole (the so-called ``ageing'' time \textit{$t_{a}$}) for
a fixed observation time $t$. The main results obtained in these
experiments are as follows. First, in almost all observation
times, the spectral hole has the shape of a Gaussian curve.
Second, the native proteins exhibit a power-law broadening of the
spectral hole by the law $\sigma_{\nu}(t)\sim t^{b}$ with the
characteristic exponent $b=0.27\pm0.03$. Third, the exponent $b$
is almost independent of the temperature $T$ of a sample in the
temperature interval from $0.1$ to $4.2$ K. Moreover, in the
experiments we also observed a power-law dependence
$\sigma_{\nu}(t_{a})\sim t_{a}^{-c}$ of the spectral hole width on
the ageing time of the sample with the characteristic exponent
\textit{$c=0.07\pm0.01$} (for a fixed observation time and a fixed
temperature of $4.2$ K.) Nevertheless, the results of the
experiments do not allow us to make a conclusion about the
dependence of the exponent $c$ on temperature.

To describe the dependence of the characteristic width of the
spectral hole on the observation time $t$, the authors of the
experimental studies \cite{LW,PFVBB,Friedrich1,Friedrich2}
proposed a model of ``superposition'' of two random processes.
This model is based on a Wiener random process on the absorption
frequency axis in which the role of time is played by the maximum
deviation of the coordinate of another Wiener random process. In
spite of the fact that exact mathematical implementation of such a
process was not presented in these works, the authors showed that
the root-mean-square deviation of the absorption frequency in such
a process should depend on time as $\sigma(t)\sim t^{0.25}$. Note
also that this model carries no information about the dynamics of
conformational rearrangements of the surrounding of a chromophore
marker and does not display the characteristic features of the
ultrametric kinetics of protein. The dependence of the
characteristic width of the spectral hole on ageing time is not
described in this model either. In \cite{AB_2008}, the authors
attempted to describe the dependence of the characteristic width
of the spectral hole on the observation time and ageing time using
the ultrametric concept of conformational dynamics of protein. For
a fixed value of the temperature parameter $\alpha$, they could
correlate the power law of $\sigma_{\nu}$ as a function of the
observation time $t$ and ageing time $t_{a}$ with experiment.
Nevertheless, this model is not completely consistent because,
first, it contains a number of hard-to-interpret assumptions and,
second, it cannot explain the independence of the exponent $b$ in
the broadening law $\sigma_{\nu}(t)\sim t^{b}$ on temperature.
Thus, at present there is no completely consistent model
describing experiments on spectral diffusion; this fact stimulates
further investigations in this direction.

Here we make an attempt to present another model for a possible
description of experiments of this kind. The main assumption of
the model is that the conformational dynamics of protein is
determined by a random walk on the ultrametric conformational
space of quasiequilibrium macrostates, which is identified with
the field of $p$-adic numbers $\mathbb{Q}_{p}$. The distribution
function of proteins over the space of conformational states is
subject to equation (\ref{UDE}), where $\alpha=\dfrac{T_{0}}{T}$
is a temperature parameter, $T$ is temperature, $T_{0}$ is a
parameter describing a temperature scale. In this case, the
variation of the absorption frequency of the chromophore marker
occurs if and only if the protein is in the conformational
subspace described by a ring of $p$-adic integers
$\mathbb{Z}_{p}$. When the protein is outside the conformational
domain $\mathbb{Z}_{p}$, no variation of the absorption frequency
of the chromophore marker occurs. To implement such a model, one
can describe the variation of the absorption frequency of the
chromophore marker by some non-Markovian Gaussian random process
$\nu\left(\theta\right)$ ($\theta$ is the time variable) on the
frequency axis. We impose the self-similarity condition
$f\left(\nu,t\right)=t^{-h}f\left(\nu t^{-h},1\right)$ with
$h\neq1$ on the distribution function of such a process
$f\left(\nu,t\right)$ and require the finiteness of the second
moment (see, for example, \cite{LM,Uch}). It follows from these
conditions that the variance of such a process depends on time
$\theta$ as $\left\langle \nu^{2}\left(\theta\right)\right\rangle
_{\nu}=D\theta^{2h}$ for some parameter $D$. Below we assume that
this process is possible if and only if the protein states belong
to the conformational domain $\mathbb{Z}_{p}$. Thus, the
dependence of the frequency $\nu\left(t\right)$ on the real time
$t$ is described by the superposition of two random processes: the
process $\nu\left(\theta\right)$ and the process
$\theta\left(t\right)$ -- the sojourn time of protein in the
conformational domain $\mathbb{Z}_{p}$ , i.e.,
$\nu\left(t\right)\equiv\nu\left(\theta\left(t\right)\right).$ In
this case the width of the spectral hole $\sigma\left(t\right)$ is
the root-mean-square deviation of the trajectory of the process
from its initial position, i.e.,
$\sigma\left(t\right)=\left(\left\langle \left\langle
\nu^{2}\left(\theta\left(t\right)\right)\right\rangle
_{\nu}\right\rangle _{\theta}\right)^{\tfrac{1}{2}}$
$=\left(\left\langle D\theta^{2h}\left(t\right)\right\rangle
_{\theta}\right)^{\tfrac{1}{2}}$. For large times $t$ and low
temperatures ($\alpha>1$), according to Theorem 6, we have
$\sigma\left(t\right)\sim
t^{\tfrac{\left(\alpha-1\right)h}{\alpha}}$. In the limit as
$\alpha\rightarrow\infty$ (at ultra-low temperatures), the
exponent $\sigma\left(t\right)$ does not depend on $\alpha$ and is
equal to $h$, which gives agreement with experiment for $h=b$.
Note that this model automatically guarantees the Gaussian shape
of the spectral hole. In addition, this makes it possible to trace
the ``fine'' temperature dependence of the exponent in the
experimentally established variation law of the spectral hole
width. In the case of ageing (i.e., when between the preparation
time of the sample and burnout time of the spectral hole $t_{a}$)
passes, the sojourn time of the protein in the conformational
subspace $\mathbb{Z}_{p}$ starting from $t_{a}$ is
$\theta\left(t_{a}+t\right)-\theta\left(t_{a}\right)$. In this
case $\sigma\left(t,t_{a}\right)=\left(\left\langle \left\langle
\nu^{2}\left(\theta\left(t_{a}+t\right)-\theta\left(t_{a}\right)\right)\right\rangle
_{\nu}\right\rangle _{\theta}\right)^{\tfrac{1}{2}}$
$=\left(D\left\langle
\left(\theta\left(t_{a}+t\right)-\theta\left(t_{a}\right)\right)^{2h}\right\rangle
_{\theta}\right)^{\tfrac{1}{2}}$. Our preliminary estimates of
this expression (the details of which we omit here) show that
$\sigma\left(t,t_{a}\right)\sim t_{a}^{-\tfrac{h}{\alpha}}$ for
$t_{a}\gg t$. In this case, the coincidence with the experimental
value of the exponent \textit{$c$} occurs at the values
$\alpha=\dfrac{b}{c}$ of the temperature parameter, which allows
one to establish a correlation between $\alpha$ and $T$ and thus
determine the value of the parameter $T_{0}$ which defines the
temperature scale. In conclusion, note that a more detailed
analysis of the model is of undoubted interest, which, in turn,
requires its precise mathematical implementation. We hope to
implement this investigation in a separate publication in the
nearest future.

\vspace {5mm}
{\bf Acknowledgments}

The study was supported in part by the Ministry of Education and Science
of Russia by State assignment to educational and research institutions
under project FSSS-2020-0014.


\begin{thebibliography}{10}
\bibitem{Takacs_1} L. Takacs, ``On certain sojourn time problems
in the theory of stochastic processes'', Acta Mathematica Academiae
Scientiarum Hungarica 8(1-2), 169-191 (1957). \href{ttps://doi.org/10.1007/BF02025241}{ttps://doi.org/10.1007/BF02025241}

\bibitem{Takacs_2} L. Takacs, ``On a sojourn time problem'', Theory
of Probability \& Its Applications 3(1), 58--65, (1958). \href{https://doi.org/10.1137/1103003}{https://doi.org/10.1137/1103003}

\bibitem{Takacs_3} L. Takacs, ``Sojourn time problems'', The Annals
of Probability, 420--431, (1974). \href{\%20https://doi.org/10.1214/aop/1176996657}{ https://doi.org/10.1214/aop/1176996657}

\bibitem{BP} R. E. Barlow, \& F. Proschan, \textit{Mathematical theory
of reliability} (Society for Industrial and Applied Mathematics, 1996).

\bibitem{Io} M. Iosifescu, \textit{Finite Markov processes and their
applications} (Courier Corporation, 2014).

\bibitem{Gr} F. Grabski, \textit{Semi-Markov processes: applications
in system reliability and maintenance}, Vol. 599 ( Amsterdam, The
Netherlands: Elsevier, 2015).

\bibitem{RTV} R. Rammal, G. Toulouse, M.A. Virasoro, ``Ultrametricity
for physicists'', Rev. Mod. Phys. 58(3), 765--788, (1986). \href{https://doi.org/10.1103/RevModPhys.58.765}{https://doi.org/10.1103/RevModPhys.58.765}

\bibitem{VVZ}  V. S. Vladimirov, I. V. Volovich, E.I. Zelenov, \textit{$p$-Adic
analysis and mathematical physics} (World Sci. Publishing, Singapore,
1994).

\bibitem{ALL}  B. Dragovich, A.Yu. Khrennikov, S.V. Kozyrev, I.V.
Volovich, ``On $p$-adic mathematical physics'', $p$-Adic Numbers,
Ultrametric Analysis and Applications 1(1), 1--17, (2009). \href{https://doi.org/10.1134/S2070046609010014}{https://doi.org/10.1134/S2070046609010014}.

\bibitem{ALL_1} B. Dragovich, A.Yu. Khrennikov, S.V. Kozyrev, I.V.
Volovich, E.I. Zelenov, ``$p$-Adic mathematical physics: The first
30 years'', $p$-Adic Numbers, Ultrametric Analysis and
Applications 9(2), 87--121, (2017).
\href{https://doi.org/10.1134/S2070046617020017}{https://doi.org/10.1134/S2070046617020017}

\bibitem{DKKM}  B. Dragovich, A.Y. Khrennikov, S. V. Kozyrev, N.
\v{Z}. Mi\v{s}i\'{c}, ``$p$-Adic mathematics and theoretical biology''.
Biosystems 199, 104288--104288, (2021). \href{https://doi.org/10.1016/j.biosystems.2020.104288}{https://doi.org/10.1016/j.biosystems.2020.104288}

\bibitem{ABK_1999} V. A. Avetisov, A. Kh. Bikulov, S. V. Kozyrev,
``Application of $p$-adic analysis to models of spontaneous breaking
of replica symmetry'', Journal of Physics A: Mathematical and General
32(50), (1999) 8785--8791. \href{https://doi.org/10.1088/0305-4470/32/50/301}{https://doi.org/10.1088/0305-4470/32/50/301}.

\bibitem{ABKO_2002} V. A. Avetisov, A. Kh. Bikulov, S. V. Kozyrev,
V. A. Osipov, ``$p$-Adic Models of ultrametric diffusion constrained
by hierarchical energy landscapes'', Journal of Physics A: Mathematical
and General 35(2), (2002), 177--189. \href{https://doi.org/10.1088/0305-4470/35/2/301}{https://doi.org/10.1088/0305-4470/35/2/301}.

\bibitem{ABO_2003} V. A. Avetisov, A. Kh. Bikulov, V. A. Osipov,
``$p$-Adic description of characteristic relaxation in complex systems'',
Journal of Physics A: Mathematical and General 36(15), 4239--4246?
(2003). \href{https://doi.org/10.1088/0305-4470/36/15/301}{https://doi.org/10.1088/0305-4470/36/15/301}.

\bibitem{ABO_2004} V. A. Avetisov, A. Kh. Bikulov, V. A. Osipov,
``$p$-Adic Models for ultrametric diffusion in conformational dynamics
of macromolecules'', Tr. Mat. Inst. Steklova 245, 55--64, (2004).
\href{\%5C\%5C\%5C\%5C\%5C\%5C\%5C\%5C\%5C\%5C\%5C\%5C\%5C\%5C\%5C\%20http://mi.mathnet.ru/eng/tm/v245/p55}{ http://mi.mathnet.ru/eng/tm/v245/p55}.

\bibitem{AB_2008} V. A. Avetisov, A. Kh. Bikulov, ``Protein ultrametricity
and spectral diffusion in deeply frozen proteins'', Biophysical Reviews
and Letters 3, 387--396, (2008). \href{https://doi.org/10.1142/S1793048008000836}{https://doi.org/10.1142/S1793048008000836}.

\bibitem{ABZ_2009} V. A. Avetisov, A. Kh. Bikulov, A. P. Zubarev,
``First passage time distribution and number of returns for ultrametric
random walk'', Journal of Physics A: Mathematical and General 42(8),
85005--85021, (2009). \href{https://doi.org/10.1088/1751-8113/42/8/085003}{https://doi.org/10.1088/1751-8113/42/8/085003}.

\bibitem{ABZ_2011} V.A. Avetisov, A.Kh. Bikulov, A.P. Zubarev, Mathematical
modeling of molecular ``nano-machines'', Vestn. Samar. Gos. Tekhn.
Univ. Ser. Fiz.-Mat. Nauki 1(22) (2011) 9-15. \href{https://doi.org/10.14498/vsgtu906}{https://doi.org/10.14498/vsgtu906}.

\bibitem{ABZ_2013} V. A. Avetisov, A. Kh. Bikulov, A. P. Zubarev,
``Ultrametricity as a basis for organization of protein molecules:
CO binding to myoglobin'', Vestn. Samar. Gos. Tekhn. Univ. Ser. Fiz.-Mat.
Nauki, 1(30), 315--325, (2013). \href{https://doi.org/10.14498/vsgtu1154}{https://doi.org/10.14498/vsgtu1154}.

\bibitem{BZ_2021} A. Kh. Bikulov, A. P. Zubarev, ``Ultrametric theory
of conformational dynamics of protein molecules in a functional state
and the description of experiments on the kinetics of CO binding to
myoglobin'', Physica A: Statistical Mechanics and its Applications,
583, 126280--126280, (2021). \href{https://doi.org/10.1016/j.physa.2021.126280}{https://doi.org/10.1016/j.physa.2021.126280}

\bibitem{Takacs_4} L. Takacs, ``On limiting distributions concerning
a sojourn time problem'', Acta Mathematica Academiae Scientiarum
Hungarica 8(3-4), 279--294, (1957). \href{https://doi.org/10.1007/bf02020316}{https://doi.org/10.1007/bf02020316}

\bibitem{LW} D. Th. Leeson, D. A. Wiersma, ?Looking into the energy
landscape off myoglobin?, Nature Struct. Biol., 2, 848--851, (1995).
\href{https://doi.org/10.1038/nsb1095-848}{https://doi.org/10.1038/nsb1095-848}

\bibitem{PFVBB} V. V. Ponkratov, J. Friedrich, K.M. Vanderkooi, A.
L. Burin, Yu. A. Berlin, ?Physics of protein at low temperature?,
J. Low. Temp. Phys. 3, 289--317, (2006). \href{https://doi.org/10.1023/B:JOLT.0000049058.81275.72}{https://doi.org/10.1023/B:JOLT.0000049058.81275.72}

\bibitem{Friedrich1} J. Schlichter, J. Friedrich, L. Herenyi, J.
Fidy, ``Protein dynamics at low temperatures'', The Journal of Chemical
Physics 112(6), 3045--3050, (2000). \href{\%20https://doi.org/10.1063/1.480879}{ https://doi.org/10.1063/1.480879}

\bibitem{Friedrich2} J. Schlichter, K.-D. Fritsch, J. Friedrich,
J. M. Vanderkooi, ``Conformational dynamics of a low temperature
protein: Free base cytochrome-c'', The Journal of Chemical Physics
110, 3229--3234 (1999). \href{http://dx.doi.org/10.1063/1.477845}{http://dx.doi.org/10.1063/1.477845}

\bibitem{LM} S. C. Lim, S. V. Muniandy, ``Self-similar Gaussian
processes for modeling anomalous diffusion'', Physical Review E 66(2),
021114.1--021114.14, (2002). \href{https://doi.org/10.1103/PhysRevE.66.021114}{https://doi.org/10.1103/PhysRevE.66.021114}

\bibitem{Uch} V. V. Uchaikin, ``Self-similar anomalous diffusion
and Levy-stable laws'', Phys. Usp. 46 821--849, (2003). \href{http://dx.doi.org/10.1070/PU2003v046n08ABEH001324}{http://dx.doi.org/10.1070/PU2003v046n08ABEH001324}
\end{thebibliography}
\end{document}